\documentclass[reqno,a4paper]{amsart}

\usepackage{amssymb}

\newcommand\idd{\mathop {\fam 0 id}\nolimits}

\newcommand\Coeff{\mathop {\fam 0 Coeff}}
\newcommand\Hom{\mathop {\fam 0 Hom}\nolimits}

\newcommand{\oo}[1]{\mathrel{{\circ }_{#1}} }

\newcommand\Ress{\mathop{ \fam 0 Res}}

\newtheorem{thm}{Theorem}
\newtheorem{lem}[thm]{Lemma}
\newtheorem{prop}[thm]{Proposition}

\theoremstyle{definition}
\newtheorem{defn}[thm]{Definition}

\newtheorem{rem}[thm]{Remark}

\numberwithin{equation}{section}
\numberwithin{thm}{section}

\begin{document}

% \baselineskip=2\baselineskip

\title{Identities of Conformal Algebras and Pseudoalgebras}

\thanks{Partially supported by RFBR and SSc-2269.203}

\author{Pavel Kolesnikov}
\email{pavelsk@kias.re.kr}

\address{Korea Institute for Advanced Study, 207-43 Cheongnyangni 2-dong,
Dongdaemun-gu, 130-722 Seoul, Korea;}

\address{Sobolev Institute of Mathematics, Acad. Koptyug prosp. 4,
630090 Novosibirsk, Russia}

\begin{abstract}
For a given conformal algebra $C$, we write down the correspondence between identities
of the coefficient algebra $\Coeff C$
and identities of $C$ itself as of pseudoalgebra.
In particular, we write down the defining relations of Jordan,
alternative and Mal'cev conformal algebras,
and show that the analogue of Artin's Theorem does not hold
for alternative conformal algebras.
\end{abstract}

\maketitle

\section{Conformal algebras}

In this note, we present a proof of a technical statement
which concerns the relation between
identities of a conformal algebra $C$ and its coefficient algebra $\Coeff C$.
This relation was mentioned in \cite{Ro}, where some particular cases
(associativity, commutativity, Jacobi identity) were considered.
Although the approach of \cite{Ro} is quite general,
it is still technically difficult
to write down the conformal identities corresponding
to a given variety of ordinary algebras.

We propose another approach which uses the language of
pseudoproduct \cite{BDK}, in order to obtain the correspondence
between identities of $C$ and $\Coeff C$ in a very explicit form.
This approach was mentioned in \cite{BDK}, where the most
important cases (associativity, commutativity, Jacobi identity)
were considered. We prove the general statement for any
homogeneous multilinear identity. As an application, we write down
the identities of Jordan, alternative and Mal'cev conformal
algebras and derive their elementary properties.

\begin{defn}[\cite{K1}]\label{defn1}
Let $\Bbbk $ be a field of zero characteristic,
and let $\Bbbk [D]$ be the polynomial algebra in one variable.
A {\em conformal algebra\/}
$C$ is a unital left $\Bbbk[D]$-module
endowed with a family of
$\Bbbk$-bilinear operations
$(\cdot \oo{n}\cdot )$
($n$ ranges the set of non-negative integers)
satisfying the following properties:
\begin{gather}
a\oo{n} b = 0 \quad \mbox{for\ sufficiently\ large}\ n,\quad a,b\in C;
                                                    \label{P1.1} \\
Da\oo{n}b= -n a\oo{n-1} b,
  \quad a\oo{n}Db=D(a\oo{n} b)+n a\oo{n-1} b, \quad n\ge 0.
                                                \label{P1.2}
\end{gather}
The conditions (\ref{P1.1}) and (\ref{P1.2})
are called {\em locality\/} and
{\em sesqui-linearity}, respectively.
\end{defn}

This definition is a formalization of the following
structure (appeared in mathematical physics)
\cite{K2,K3}.
Let $A$ be an algebra over $\Bbbk$ (in general, $A$ is non-associative),
and let $A[[z,z^{-1}]]$
be the space of formal distributions over
$A$: $A[[z,z^{-1}]] = A \otimes \Bbbk [[z,z^{-1}]]$.
An ordered pair of distributions
$\langle a(z), b(z) \rangle $
is said to be {\em local},
if
\begin{equation}\label{locality}
a(w)b(z) (w-z)^N = 0
\end{equation}
for some $N\ge 0$.
If $\langle a(z), b(z)\rangle$ is a local pair,
then the product $a(w)b(z)\in A[[z,z^{-1},w,w^{-1}]]$
could be presented as a finite sum \cite{K1}
\begin{equation}\label{OPE}
a(w)b(z)  = \sum\limits_{n\ge 0} c_n(z) \frac{1}{n!} \partial^n_z \delta(w-z),
\end{equation}
where $\delta(w-z) = \sum_{s\in {\mathbb Z}} w^{s} z^{-s-1}$
is the formal delta-function.
In physics, the relation (\ref{OPE})
 is known as the operator product
expansion (OPE) of conformal fields.

Denote $c_n(z)$ by $(a\oo{n} b)(z)$.
The explicit expression for $c_n$
could be easily derived \cite{K1}:
\begin{equation}\label{Rez}
(a\oo{n} b)(z) = \Ress\nolimits_{w} a(w)b(z) (w-z)^n,
\end{equation}
where $\Ress_w$ means the coefficient at $w^{-1}$.
The following statement is straightforward.

\begin{thm}[see, e.g., \cite{K3, Ro}]\label{thmForDist}
Let $C$ be a subspace of  $A[[z,z^{-1}]]$
such that
\begin{itemize}
\item
any $a(z),b(z)\in C$ form a local pair;
\item
for any $a(z),b(z)\in C$ and for any $n\ge 0$ we have $(a\oo{n} b)(z) \in C$;
\item
for any $a(z)\in C$ its derivative $\frac{d}{dz}a(z)$ lies in $C$.
\end{itemize}
Then $C$ is a conformal algebra with respect to the operations
(\ref{Rez}) and $D = \frac{d}{dz}$.
\end{thm}

\begin{rem}[Dong's Lemma, \cite{K1}]\label{lem Dong}
If $A$ is an associative or Lie algebra, then
any subset $B\subset A[[z,z^{-1}]]$ such that $B\times B$
consists of local pairs generates a conformal algebra.
\end{rem}

The converse is also true: for any conformal algebra $C$
in the sense of Definition~\ref{defn1}
one can build an algebra $A$ such that
$C$ could be represented as a subspace
of $A[[z,z^{-1}]]$ satisfying the conditions of
 Theorem \ref{thmForDist}.
There exists a universal algebra of this kind
(it is unique up to isomorphism) called the {\em coefficient algebra\/}
$\Coeff C$~\cite{Ro}.

Let us remind the construction of $\Coeff C$
for a given conformal algebra~$C$
(see \cite{K3} or \cite{Ro} for details).
By the definition, $C$ is a (unital) left module over $H=\Bbbk [D]$.
Consider the vector space
 $\Bbbk [t,t^{-1}]$ as a right $H$-module with respect
to the action
$t^n D = -n t^{n-1}$.
The underlying vector space of $\Coeff C$  is
$\Bbbk[t,t^{-1}]\otimes_H C$.
Denote $t^n\otimes _H a$ by $a(n)$, $a\in C$, $n\in {\mathbb Z}$.
Define a multiplication on this space via the formula
\begin{equation}
a(n)b(m) = \sum\limits_{s\ge  0} \binom{n}{s} (a\oo{s} b)(n+m-s).
                            \label{CoeffMul}
\end{equation}
This is a well-defined bilinear operation which makes $\Coeff C$
to be an algebra.

\begin{thm}[\cite{Ro}]\label{thmCoeff}
Let $C$ be a conformal algebra, and let $A=\Coeff C $.
Then

{\rm (i)} $C$ lies in $A[[z,z^{-1}]]$ as a subspace satisfying the
conditions of Theorem \ref{thmForDist}: the embedding is given by
$a\mapsto \sum_{n\in {\mathbb Z}}a(n)z^{-n-1}$;

{\rm (ii)} for any algebra $B$ and for any homomorphism
of conformal algebras $\varphi : C\to C'\subset B[[z,z^{-1}]]$
there exists a homomorphism of algebras $\psi : A\to B$
such that $\varphi(a) = \sum_{n\in {\mathbb Z}}\psi(a(n))z^{-n-1}$,
$a\in C$.
\end{thm}

\begin{defn}[\cite{Ro}]\label{ConfVar}
Let $\Omega $ be a variety of algebras. A conformal algebra $C$
is said to be an $\Omega$-conformal algebra if and only if
$\Coeff C $ belongs to $\Omega $.
\end{defn}

It was shown in \cite{Ro}
how to convert an identity of a coefficient algebra $\Coeff C$
into the corresponding series
of identities of~$C$.
In this way, one should proceed with a routine
computation, and it is difficult to predict the final result.
The most common identities
have the following form
\cite{K3,Ro}:

\medskip
\noindent
{\bf Associativity}
\begin{equation}
(a\oo{n} b)\oo{m} c=
\sum\limits_{t\ge 0} (-1)^t \binom{n}{t} a \oo{n-t} (b \oo{m+t} c);
                                              \label{ConfAss}
\end{equation}
{\bf (Anti-)Commutativity}
\begin{equation}
 a\oo{n} b \pm
 \sum\limits_{s\ge 0}
 \frac{(-1)^{n+s}}{s!} D^{s}(b\oo{n+s} a) =0;
                                                      \label{ConfComm}
\end{equation}
{\bf Jacobi identity}
\begin{equation}
 (a\oo{n} b)\oo{m} c
 =\sum\limits_{s\ge 0}(-1)^s \binom{n}{s}
(a\oo{n-s}( b\oo{m+s} c) - b\oo{m+s} (a \oo{n-s} c)).
                                                       \label{ConfJac}
\end{equation}

\section{Pseudoalgebras}

Theorems \ref{thmForDist} and \ref{thmCoeff} show that
any conformal algebra could be considered as an algebraic
structure on a formal distribution space.
But there is a more formal approach to the theory of conformal algebras
related with the notion of a pseudo-tensor category \cite{BD}.

Consider $H=\Bbbk [D]$ as a Hopf algebra with
respect to the usual
coproduct $\Delta (D) = 1\otimes D + D\otimes 1$,
counit $\varepsilon(D)=0$, and antipode $S(D)=-D$.
We will use the standard notation
\begin{eqnarray}
& \Delta(f) =\sum\limits_i \frac{D^{i}}{i!}\otimes f^{(i)}
  = f_{(1)}\otimes f_{(2)},
  \nonumber \\
& (\Delta\otimes \idd) \Delta(f) = (\idd\otimes \Delta) \Delta(f)
    = f_{(1)}\otimes f_{(2)}\otimes f_{(3)},
 \quad \mbox{etc.}
 \nonumber
\end{eqnarray}
Let us denote by $\Delta^{(n)}$, $n\ge 1$, the iterated
coproduct: $\Delta^{(1)}=\idd_H$,
$\Delta^{(n+1)} = (\idd_H \otimes\Delta^{(n)})\Delta$.

The algebra $H$ acts on its tensor power
 $H^{\otimes n}$
as follows:
\[
(f_1\otimes \dots \otimes f_n)h = f_1h_{(1)}\otimes \dots \otimes f_nh_{(n)},
\quad f_i,h\in H.
\]

For a given conformal algebra $C$, define the operation
\begin{gather}
*: C\otimes C \to (H\otimes H)\otimes_H C,
\nonumber \\
a* b = \sum\limits_{s\ge 0} \frac{(-D)^s}{s!} \otimes 1 \otimes_ H (a\oo{s} b),
\quad a,b\in C,
\label{pseudoprod}
\end{gather}
called {\em pseudoproduct}
(the axiom (\ref{P1.1}) implies this sum to be finite).
It follows from (\ref{P1.2}) that the pseudoproduct is $H$-bilinear:
\begin{equation}\label{H-bilin}
    f(D)a*g(D)b = (f(D)\otimes g(D) \otimes _H 1)(a*b).
\end{equation}

\begin{defn}[\cite{BDK}]
Let $H$ be a bialgebra.
A left unital $H$-module $P$
endowed with an $H$-bilinear map
$*: P\otimes P \to (H\otimes H)\otimes _H P$
is called an $H$-{\em pseudoalgebra}.
\end{defn}

\begin{rem}[\cite{BDK}]
If $H=\Bbbk $, then $H$-pseudoalgebra
is just an ordinary algebra over the field~$\Bbbk $.
If $H=\Bbbk[D]$,
then we obtain the notion of a conformal algebra.
\end{rem}

One of the main features of a pseudo-tensor category is the composition
of multi-operations \cite{BD}. Any $H$-pseudoalgebra is just an algebra
 in the pseudo-tensor category
${\mathcal M}(H)$  (see , e.g., \cite{BDK} for details)
associated with~$H$.
In this category, an arbitrary composition of $H$-bilinear maps could
be described by the following structure.

\begin{defn}\label{expanded}
Let $P$ be an $H$-pseudoalgebra with a pseudoproduct~$*$.
The {\em expanded pseudoproduct\/}
is an
 $H^{\otimes (n+m)}$-linear map
\[
* : (H^{\otimes n}\otimes_H P)\otimes (H^{\otimes m}\otimes_H P) \to
 H^{\otimes (n+m)}\otimes_H P,
 \quad n,m\ge 1,
\]
defined as
\begin{equation}
(F\otimes_H a)*(G\otimes _H b) = (F\otimes G\otimes_H 1)
(\Delta^{(n)}\otimes \Delta^{(m)}\otimes_H\idd_P)(a*b),
                        \label{J1.12}
\end{equation}
$F\in H^{\otimes n}$, $G\in H^{\otimes m}$, $a,b\in P$.
\end{defn}

Now, fix $H=\Bbbk[D]$ (an $H$-pseudoalgebra is the same as
conformal algebra). The series of identities (\ref{ConfAss}),
(\ref{ConfComm}), (\ref{ConfJac}) could be expressed in terms of
the expanded pseudoproduct (\ref{J1.12}) as follows \cite{BDK}:

\medskip
\noindent
{\bf Associativity}
\begin{equation}
(a* b)*c= a*(b*c);
                                                     \label{PsAss}
\end{equation}
{\bf (Anti-)Commutativity}
\begin{equation}
 a* b \pm (\tau_{12}\otimes_H \idd_P )(b*a)=0;
                                                      \label{PsComm}
\end{equation}
{\bf Jacobi identity}
\begin{equation}
 (a* b)*c=a*( b*c) - (\tau_{12}\otimes_H \idd_P)(b*(a* c)).
                                                \label{PsJac}
\end{equation}
Here $\tau_{12}$ means the permutation of two tensor factors
in $H\otimes H $ or $H\otimes H\otimes H$ (it is a well-defined
$H$-module automorphism provided that $H$ is cocommutative).

It is easy to see that the expressions
(\ref{PsAss}), (\ref{PsComm}), (\ref{PsJac})
are similar in some sense
to the ordinary associativity, (anti-)commutativity, and Jacobi identity,
respectively. It is natural to suppose that the similarity
holds for an arbitrary identity.
In the next section, we prove the correspondence.

\section{Identities of conformal algebras}

Let us consider a conformal algebra $C$ and its coefficient algebra
$\Coeff C = \Bbbk[t,t^{-1}]\otimes _H C$ with the multiplication (\ref{CoeffMul}).
In this section, we denote $x\otimes _H a\in \Coeff C$ by $a(x)$,
$x\in \Bbbk [t,t^{-1}]$, and
$H=\Bbbk[D]$ as before. Note that the dual algebra $X=H^*$ is isomorphic
to $\Bbbk[[t]]$, where $\langle t, D^n \rangle = \delta_{n,1}$.

Consider the associative topological algebra
 $X^{(1)}= \Bbbk[t,t^{-1}] \otimes \Bbbk[t] $,
where the basic neighborhoods of zero are of the form
 $\{(1\otimes t^n)f \mid f\in X^{(1)}\}$,
 $n\ge 0$.
Denote by
$\overline{X^{(1)}}$
the completion  of $X^{(1)}$.
The algebra
$\overline{X^{(1)}}$
consists of series like
 $\sum\limits_{s\ge 0} f_s(t) \otimes t^s $,
 $f_s(t)\in \Bbbk[t,t^{-1}]$.

The standard coproduct
$\Delta : \Bbbk[t]\to \Bbbk[t]^{\otimes 2}$,
$\Delta(t)=t\otimes 1 + 1\otimes t$,
could be continued  to the homomorphism
 $\Delta :  \Bbbk[t,t^{-1}] \to \overline{X^{(1)}}$
 via
\begin{equation}
\Delta (t^n)= \sum\limits_{s\ge 0} \binom{n}{s} t^{n-s}\otimes t^{s},
\quad n\in {\mathbb Z}.
                                                              \label{LDelta}
\end{equation}
Let us denote
 $\Delta(x)=x_{(1)}\otimes x_{(2)}$, as before.

By the same way, we can define the natural topology on $X^{(r)} =
X^{(r-1)}\otimes  \Bbbk[t]$, $r>1$. Namely, let the basic
neighborhoods of zero be of the form
\[
\{(1\otimes t^{n_1}\otimes\dots \otimes t^{n_r}) f     \mid
   f\in X^{(r)},\,
 n_1+\dots +n_r = n \},
 \quad
 n\ge 0.
\]
By  $\overline{X^{(r)}}$
we denote the completion of $X^{(r)}$.
Then the map
 $\idd\otimes \Delta : X^{(1)} \to X^{(2)}$
is a continuous homomorphism, as well as
$\Delta\otimes \idd : X^{(1)}\to \overline{X^{(2)}}$.
Hence,
$\Delta \otimes \idd$
and
$\idd\otimes \Delta$
are defined on
 $\overline{X^{(1)}}$, it is easy to check that they coincide on
 $\Delta(\Bbbk[t,t^{-1}])$.
So  the ``expanded"  homomorphism
 $\Delta$ is coassociative.
 By $S$ we denote the standard antipode of $\Bbbk [t]$ given by $S(t)=-t$.

It is easy to see that (\ref{CoeffMul}) is  equivalent to
\begin{equation}
a(x)b(y) = (a \oo{x_{(2)}} b) (x_{(1)}y),
\quad x,y\in \Bbbk [t,t^{-1}],
                                                 \label{J3.5}
\end{equation}
so
\begin{equation}
(a \oo{x} b)(y) = a(x_{(2)})b(S(x_{(1)}) y),
\quad
x\in \Bbbk[t],\  y\in \Bbbk[t,t^{-1}],
                                    \label{J3.5-1}
\end{equation}
for any
$a,b\in C$.

Let $f(a_1,\dots, a_n)$ be a (non-associative) homogeneous
multilinear polynomial with coefficients in $\Bbbk $. Any
polynomial of this kind could be written as
\begin{equation}\label{J3.8}
f(a_1,\dots, a_n) = \sum\limits_{\sigma \in {\mathbb S}_n}
t_\sigma (a_{1\sigma},\dots ,a_{n\sigma}),
\end{equation}
where each of the terms
$t_\sigma(b_1,\dots ,b_n)$
is a linear combination of non-associative words obtained
from  $b_1\dots b_n$
by some bracketings.
For any (non-associative) monomial $m(b_1,\dots ,b_n)$ in (\ref{J3.8}),
replace the usual multiplication with the (expanded) pseudoproduct~$*$.
We obtain an expression $m^*(b_1,\dots, b_n)$
which has sense in a pseudoalgebra.
Denote by $f^*(a_1,\dots, a_n)$ the result of this operation:
\begin{equation}
f^*(a_1,\dots , a_n) =
\sum\limits_{\sigma \in {\mathbb S}_n}  (\sigma \otimes_H \idd_C)
t^*_\sigma (a_{1\sigma},\dots ,a_{n\sigma}).
                              \label{J3.8'}
\end{equation}

\begin{thm}\label{thm 3.1}
Let $C$ be a conformal algebra such that $\Coeff C$ satisfies
homogeneous multilinear identity $f=0$. Then $C$ as a
pseudoalgebra satisfies the  (pseudo-) identity
 $f^*=0$.
\end{thm}

\begin{proof}
Let us denote by
 $\theta $ the following map
 $H^{\otimes n}\otimes_H C \to H^{\otimes (n-1)}\otimes C$:
\[
\theta((h_1\otimes \dots \otimes h_{n-1} \otimes 1)\otimes_H c) =
h_1\otimes \dots \otimes h_{n-1} \otimes c.
\]
For any $x_i \in \Bbbk[t]\subset X = \Bbbk [[t]]$,
$i=1,\dots, n-1$, $n\ge 1$,
define the linear operator
\begin{equation}
 P^{x_1,\dots ,x_{n-1}} = (\langle S(x_1),\cdot \rangle \otimes
\dots \otimes \langle S(x_{n-1}),\cdot \rangle \otimes \idd_C) \theta
:  H^{\otimes n}\otimes _H C \to C,
                                                  \label{J3.7}
\end{equation}
Therefore,
\[
P : \Bbbk[t]^{n-1} \to \Hom (H^{\otimes n}\otimes _H C, C)
\]
provided by $(x_1,\dots, x_{n-1}) \mapsto P^{x_1,\dots, x_{n-1}}$
is a multilinear map which could be defined on $\Bbbk[t]^{\otimes
n-1}$. Since $P$ is continuous (with respect to the finite
topology on $\Hom (H^{\otimes n}\otimes _H C, C)$),
one may build a map
\[
  \overline{X^{(n-1)}} \to \Hom (H^{\otimes n}\otimes _H C, \Coeff C)
\]
which sends $(x_1\otimes x_2 \otimes \dots \otimes x_n)\in X^{(n-1)}$
to
\[
P^{x_2\otimes \dots \otimes x_n}_{x_1}:
 A \mapsto [P^{x_2\otimes \dots \otimes x_n}(A)](x_1),
\quad A\in H^{\otimes n}\otimes_H C.
\]

It is clear that an element $A\in H^{\otimes n}\otimes_H C$ is zero
if and only if
$P^{\bar x}(A) = 0$
for any
$\bar x \in \Bbbk[t]^{\otimes n-1}$.
Also, an element $a\in C$ is equal to zero if and only if
 $a(y) = 0$
 for any
 $y\in \Bbbk[t,t^{-1}]$,
 see Theorem~\ref{thmCoeff}.

We will use the following notation in order to simplify computations.
For
$\bar x\in \Bbbk[t]^{\otimes n-1}$ set
$\Pi \bar x$ to be the product of its components:
$\Pi(x_1\otimes \dots \otimes x_{n-1}) = x_1\dots x_{n-1}\in \Bbbk[t]$,
and denote by
$\Delta : \bar x \mapsto \bar x_{(1)}\otimes \bar x_{(2)}$
the componentwise coproduct.
For any $\bar a \in C^{\otimes n}$,
$\bar x \in \Bbbk[t,t^{-1}] ^{\otimes n}$
set
$\bar a (\bar x) \in (\Coeff C)^{\otimes n} $
in the obvious way.

\begin{lem}\label{lem 3.0}
Let $A\in H^{\otimes s}\otimes_H C$,
$B\in H^{\otimes n-s}\otimes_H C$,
$\bar x\in \Bbbk[t]^{\otimes s-1}$,
$\bar z\in \Bbbk [t]^{\otimes n-s-1}$,
$y\in \Bbbk[t]$.
Then
\begin{equation}
P^{\bar x \otimes y \otimes \bar z}(A*B) = P^{y \Pi \bar x_{(1)} }
\big ( P^{\bar x_{(2)}}(A) * P^{\bar z}(B) \big).
                                             \label{MulP}
\end{equation}
\end{lem}

\begin{proof}
We may assume $A=F\otimes 1 \otimes_H a$, $B=G\otimes 1\otimes _H b$,
and let $a*b =\sum_i h_i \otimes 1 \otimes _H c_i$.
Then
\[
A*B = \sum_i Fh_{i(1)}\otimes h_{i(2)}\otimes G \otimes 1 \otimes_H c_i,
\]
in accordance with (\ref{J1.12}).
The left-hand side of (\ref{MulP}) could be expressed as
\begin{multline}
\sum_i  \langle S(\bar x), Fh_{i(1)} \rangle
\langle S(y), h_{i(2)}\rangle
\langle S(\bar z), G \rangle  c_i  \\
= \sum_i  \langle S(\bar x_{(2)}), F\rangle
 \langle \Pi S(\bar x_{(1)}), h_{i(1)} \rangle
 \langle S(y), h_{i(2)}\rangle
 \langle S(\bar z), G \rangle  c_i \\
   = \sum_i  \langle S(\bar x_{(2)}), F \rangle
      \langle S(y\Pi \bar x_{(1)}), h_i \rangle
      \langle S(\bar z), G\rangle   c_i,
\end{multline}
which is equal to the right-hand side of (\ref{MulP}).
\end{proof}

\begin{lem}\label{lem 3.1}
Let $t(a_1,\dots, a_n)$ be a non-associative monomial
obtained from $a_1\dots a_n$
by some bracketing.

{\rm (i)}
For any
$\bar x \in \Bbbk[t]^{\otimes n-1}$,
 $y\in \Bbbk[t,t^{-1}]$,
 $\bar a \in C^{\otimes n-1}$,
we have
\begin{equation}
P^{\bar x}_y (t^*(\bar a, a_n))
=t(\bar a (\bar x_{(2)}), a_n(S(y\Pi \bar x_{(1)}) )).
                              \label{J3.10}
\end{equation}

\noindent
{\rm (ii)}
For any
$\bar a \in C^{\otimes n }$ and for any
$\bar x =x_1\otimes \dots \otimes x_{n}\in \Bbbk[t,t^{-1}]^{\otimes n}$
we have
\begin{equation}
t(\bar a(\bar x) )
 = P^{x_{1(2)}\otimes \dots \otimes x_{n-1(2)}}
     _{x_{1(1)}\dots x_{n-1(1)}x_n}
 t^*(\bar a)
                                                            \label{Cntr}
\end{equation}
\end{lem}

\begin{proof}
(i)
If $n=2$, then (\ref{J3.10}) coincides with (\ref{J3.5-1}).
If $n>2$ then we may assume that there is a non-trivial
decomposition
\[
t= t_1(\bar a, a_s)t_2(\bar b, b_m),
\]
where
$\bar a = (a_1,  \dots ,  a_{s-1})$,
$\bar b = (b_{1},  \dots ,  b_{m-1})$.
Let
$\bar x \in \Bbbk[t]^{\otimes s-1}$,
$x_s \in \Bbbk[t]$,
$\bar z \in \Bbbk[t]^{\otimes m-1}$,
$y\in \Bbbk[t,t^{-1}]$.
Then by Lemma \ref{lem 3.0} and by (\ref{J3.5}) we have
\begin{multline}
P^{\bar x, x_s, \bar z}_{y}
   (t_1^*(\bar a, a_s) * t_2^*(\bar b, b_m) )
 =
P^{x_s\Pi \bar x_{(1)}}_{y}
(
P^{\bar x_{(2)}}t_1(\bar a, a_s) * P^{\bar z} t_2(\bar b, b_m)
)            \\
=
\bigg(
 P^{\bar x_{(3)}}_{x_{s(2)}\Pi \bar x_{(2)}}
   t_1(\bar a, a_s) \bigg )
\bigg (
P^{\bar z}_{y S(x_{s(1)})S(\Pi \bar x_{(1)})}
 t_2(\bar b, b_m)
  \bigg ).
\end{multline}
The inductive assumption allows to proceed as follows:
\begin{multline}
P^{\bar x, x_s, \bar z}_{y}
   (t_1^*(\bar a, a_s) * t_2^*(\bar b, b_m) )
    =
  t_1(\bar a(\bar x_{(4)}), a_s(x_{s(2)}\Pi \bar x_{(2)}S(\Pi \bar x_{(3)})))   \\
\times
  t_2(\bar b (\bar z_{(2)}), b_m(y S(x_{s(1)})S(\Pi \bar x_{(1)})S(\bar z_{(1)}) )  )
   =  t_1(\bar a(\bar x_{(2)}), a_s(x_{s(2)})  )
     \\
\times
  t_2(\bar b(\bar z_{(2)}), b_m (y S(\Pi \bar x_{(1)}x_{s(1)}\Pi \bar z_{(1)}))),
\end{multline}
so we obtain~ (\ref{J3.10}).

(ii)
The proof is analogous to the one of (i).
\end{proof}

For a given $\bar x\in \Bbbk[t]^{\otimes n-1}$,
$y\in \Bbbk[t,t^{-1}]$,
and for any permutation $\sigma \in {\mathbb S}_n$,
define
$\bar \zeta(\bar x, y, \sigma)\in \Bbbk[t,t^{-1}]^{\otimes n}$
as follows:
$\bar \zeta = \zeta_1\otimes \dots \otimes \zeta_n$,
where
$\zeta_i = x_{i\sigma(2)}$     for    $i\sigma\neq n$,
$\zeta_{s} = y S(\Pi \bar x_{(1)})$,    for   $s = n\sigma^{-1}$.
The following statement generalizes Lemma \ref{lem 3.1}

\begin{lem}\label{lem 3.2}
For any $\bar x \in \Bbbk[t]^{\otimes n-1}$, $y\in \Bbbk [t,t^{-1}]$,
$\bar a \in C^{n}$, and for any $\sigma \in {\mathbb S}_n$
we have
\begin{equation}
P^{\bar x}_{y}
 ((\sigma \otimes_H \idd_C)t^*(\bar a_{\sigma }) )
=
t ( \bar a_{\sigma } (\bar \zeta(\bar x, y, \sigma))  ),
                                                   \label{J3.11}
\end{equation}
where
$\bar a_{\sigma}(\bar \zeta)
 =a_{1\sigma}(\zeta_1)\otimes \dots \otimes a_{n\sigma}(\zeta_n))$.
\end{lem}

\begin{proof}
First, let us assume
$n\sigma =n$. Then it is straightforward to check that
for any $A\in H^{\otimes n}\otimes_H C$
and for any
$\bar x=x_1\otimes \dots \otimes x_{n-1}\in \Bbbk[t]^{\otimes n-1}$
we have
\begin{equation}
P^{\bar x } ((\sigma \otimes_H \idd_C)A)
= P^{\bar x_ \sigma } A,
                                                      \label{Sigma1}
\end{equation}
where $\bar x_\sigma = x_{1\sigma}\otimes \dots \otimes x_{(n-1)\sigma}$.

Second, assume $\sigma $ to be a transposition
$\sigma =(s\, n)$, $s\in \{1,\dots ,n-1\}$.
Then
for any $A\in H^{\otimes n}\otimes_H C$,
$\bar x\in \Bbbk[t]^{\otimes s-1}$, $x_s\in \Bbbk[t]$,
$\bar z\in \Bbbk[t]^{\otimes n-s-1}$,
$y\in \Bbbk[t,t^{-1}]$,
we have
\begin{equation}
P^{\bar x\otimes x_s \otimes \bar z}_{y} ((\sigma\otimes_H \idd_C)A)
=
P^{\bar x_{(2)}\otimes S(\Pi \bar x_{(1)} \Pi\bar z_{(1)} x_s )y_{(2)}
 \otimes    \bar z_{(2)} }_{y_{(1)}}
 A.
                            \label{Sigma2}
\end{equation}
To obtain the last relation, one needs the equality
\[
\langle x, h_{(1)} \rangle (yh_{(2)}) = \langle xS(y_{(2)}), h\rangle y_{(1)},
\]
which is straightforward to check for any
$h\in H$, $x\in \Bbbk[t]$, $y\in \Bbbk[t,t^{-1}]$.

An arbitrary permutation $\sigma \in {\mathbb S}_n$, $n\sigma=s$,
could be presented as
$\sigma = (s\, n)\sigma_1$, where $n \sigma_1 = n$.
It follows from (\ref{Sigma1}), (\ref{Sigma2}) that
for $n\neq s$ we have
\begin{multline}
P^{\bar x^1 \otimes x_s \otimes \bar x^2}_{y}
 ((\sigma \otimes_H \idd_C)A )   \\
=
P^{\bar x_{(2)} \otimes
   S(\Pi \bar x^1_{(1)} \Pi\bar x^2_{(1)} x_s )y_{(2)} \otimes
   \bar x^2_{(2)} }_{y_{(1)}}
((\sigma_1 \otimes_H \idd_C)A) =
P^{\zeta' (y_{(2)}) }_{y_{(1)}} A ,
                        \label{J3.11A}
\end{multline}
where $\bar x= \bar x^1\otimes x_s\otimes \bar x^2 $,
$\bar x^1 \in \Bbbk [t]^{\otimes s-1}$,
$\bar x^2 \in \Bbbk[t]^{\otimes n-s-1}$,
and
\[
\zeta'(y) =
(\bar x_{(2)} \otimes
S(\Pi \bar x^1_{(1)} \Pi\bar x^2_{(1)} x_s )y \otimes \bar x^2_{(2)})_{\sigma_1}.
\]

Now, compute the right-hand side of (\ref{J3.11}) via (\ref{Cntr}).
For $n\sigma =n$ or $n\sigma = s\neq n$,
it coincides with the right-hand side of (\ref{Sigma1})
or  (\ref{J3.11A}), respectively, if we substitute
 $A=t^*(a_{1\sigma}, \dots, a_{n\sigma})$.
\end{proof}

Let us complete the proof of Theorem~3.1.
Lemmas \ref{lem 3.1} and \ref{lem 3.2} imply that if
$\Coeff C$
satisfies the identity
 (\ref{J3.8}),
 then
\[
P^{x_1\otimes\dots\otimes x_{n-1}} _{y}
 \sum\limits_{\sigma \in {\mathbb S}_n} (\sigma
\otimes_H\idd_C) t_{\sigma}^*(a_{1\sigma} ,\dots ,a_{n\sigma})
  = 0
\]
for any $x_1,\dots, x_{n-1}\in \Bbbk [t]$,
$y\in \Bbbk [t^{-1}, t]$.
Hence,
$C$
satisfies  (\ref{J3.8'}) as a pseudoalgebra.
\end{proof}

\begin{thm}\label{thm 3.2}
Let  $C$ be a conformal algebra.
If  $C$ as a pseudoalgebra satisfies a
(pseudo-) identity $f^* = 0$
of the form
(\ref{J3.8'}), then
$\Coeff C$
satisfies the corresponding identity $f = 0$.
\end{thm}

\begin{proof}
It is sufficient to show that
\begin{equation}
t(a_{1\sigma}(x_{1\sigma }), \dots,   a_{n\sigma}(x_{n\sigma }))
=
P^{x_{1(2)}, \dots, x_{n-1(2)} } _{ x_n x_{1(1)}\dots x_{n-1(1)} }
  (\sigma \otimes_H \idd_C)
t^*(a_{1\sigma }, \dots , a_{n\sigma} ) ,
                                                                \label{Cntr2}
\end{equation}
for any
$x_i\in \Bbbk [t,t^{-1}]$,
$\sigma \in {\mathbb S}_n$.

Indeed,
let us rewrite the left-hand side of
(\ref{Cntr2}) by (\ref{Cntr}):
in the obvious notations, we obtain
\begin{equation}
t(\bar a_\sigma (\bar x_\sigma))
 = P^{\bar y_{(2)}}_{x_s \Pi \bar y_{(1)}} t^*(\bar a_\sigma),
                                \label{Cntr3}
\end{equation}
where
$\bar x_\sigma = (\bar y \otimes x_s)$, for $s=n\sigma$.

The right-hand side of (\ref{Cntr2}) could be rewritten by
(\ref{Sigma1}) or (\ref{J3.11A}). It is easy to note that the expression obtained
coincides with (\ref{Cntr3}).

The final statement  follows directly from
(\ref{Cntr}) and (\ref{Cntr2}).
\end{proof}

We will need the identities analogous to
(\ref{J3.8'}) for the expanded pseudoproduct.
Let $n\ge 1$, and let
$\pi = \{m_i \mid i=1,\dots, n\}$
be a family of positive integers.
For a given $\sigma \in {\mathbb S}_n$,
define $ \sigma_\pi \in {\mathbb S}_{m_1+\dots + m_n}$
in such a way that
\begin{equation}
\sigma_\pi (\Delta^{(m_{1\sigma })} \otimes \dots
\otimes\Delta^{(m_{n\sigma })})(F) = (\Delta^{(m_1)} \otimes \dots
\otimes\Delta^{(m_n)}) \sigma (F)
                        \label{Vareps}
\end{equation}
for any $F\in H^{\otimes n}$.

\begin{prop}\label{prop 3.3}
Let $C$ be a conformal algebra satisfying
 a homogeneous multilinear identity
\[
\sum\limits_{\sigma \in {\mathbb S}_n} (\sigma \otimes_H \idd_C)
t_\sigma ^*(a_{1\sigma }, \dots , a_{n\sigma }) = 0.
\]
Then for the expanded pseudoproduct
 (\ref{J1.12}) we have
\begin{eqnarray}
& {\displaystyle \sum\limits_{\sigma \in {\mathbb S}_n}   }
   (\sigma_\pi \otimes_H \idd_C)
   t_\sigma ^*(A_{1\sigma }, \dots , A_{n\sigma }) = 0,
                                       \label{J3expanded} \\
& A_i = G_i\otimes_H a_i, \quad
   a_i\in C, \quad  G_i\in H^{\otimes m_i},\
   m_i \ge 1 ,\ i=1,\dots, n.
                        \nonumber
\end{eqnarray}
\end{prop}

\begin{proof}
It follows from the definition of expanded pseudoproduct
(\ref{J1.12}) that for a (non-associative) homogeneous multilinear
term $t$ we have
\[
t^*(A_1,\dots ,A_n)
= (G_1\otimes \dots \otimes G_n\otimes_H 1) (\Delta^{(m_1)}
\otimes \dots \otimes\Delta^{(m_n)} \otimes_H \idd_C) t^*(a_1,\dots ,a_n).
\]
Hence,
\begin{multline}
t^*(A_{1\sigma },\dots ,A_{n\sigma })   \\
= (G_{1\sigma }\otimes \dots \otimes G_{n\sigma }\otimes_H 1)
(\Delta^{(m_{1\sigma })} \otimes \dots \otimes\Delta^{(m_{n\sigma })} \otimes_H \idd_C)
 t^*(a_{1\sigma },\dots ,a_{n\sigma }).
\nonumber
\end{multline}
Since (\ref{Vareps}), we have
\begin{multline}
\sum\limits_{\sigma \in {\mathbb S}_n} (\sigma_\pi \otimes_H \idd_C )
   t_\sigma ^*(A_{1\sigma }, \dots , A_{n\sigma })
=
(G_1\otimes \dots \otimes G_n\otimes_H 1)            \\
\times
(\Delta^{(m_{1})} \otimes \dots \otimes \Delta^{(m_{n})}
\otimes_H \idd_C)
 \sum\limits_{\sigma \in {\mathbb S}_n} (\sigma \otimes_H \idd_C )
  t_\sigma ^*(a_{1\sigma }, \dots , a_{n\sigma }) = 0,
\end{multline}
and
(\ref{J3expanded})
holds.
\end{proof}

\begin{rem}\label{rem3x}
Note that an analogue of Theorem \ref{thm 3.2}
also holds for an arbitrary pseudoalgebra
over a cocommutative Hopf algebra $H$ \cite{BDK}:
if a pseudoalgebra $P$ satisfies an identity $f^*=0$
 of the form (\ref{J3.8'}), then
its annihilation algebra \cite{BDK} ${\mathcal A}(P)$
satisfies $f=0$.
\end{rem}

\section{Jordan, alternative and Mal'cev conformal algebras}

Theorems \ref{thm 3.1} and \ref{thm 3.2} together
with Remark~\ref{rem3x} provide a foundation for
the following generalization of Definition~\ref{ConfVar}.

\begin{defn}
Let $\Omega $ be a variety of ordinary algebras defined by a
family of homogeneous multilinear identities
$\{f_i(x_1,\dots,x_{n_i})=0\}_{i\in I}$,
 and let $P$
be a pseudoalgebra over a cocommutative Hopf algebra~$H$.
If $P$ satisfies the identities
$\big\{ f_i^* (x_1,\dots, x_{n_i})=0 \big\}_{i\in I}$,
then $P$ is said to be an $\Omega $-pseudoalgebra.
\end{defn}

In this section, we write down the identities
of Jordan, alternative and Mal'cev conformal algebras
obtained by Theorems \ref{thm 3.1} and \ref{thm 3.2}.

Let $A$ be an algebra over a field $\Bbbk$
with bilinear multiplication
$\cdot : A\times A \to A$.

\begin{defn}[see, e.g., \cite{ZSSS}]\label{defnJordan}
A commutative algebra
$(A, \cdot)$  is said to be {\em Jordan},
if it
satisfies
\begin{equation}
((a\cdot a)\cdot b)\cdot a = (a\cdot a)\cdot (b\cdot a).
                              \label{J2.5}
\end{equation}
\end{defn}

In the multilinear form (remind that  $\mbox{char}\,\Bbbk =0$),
the Jordan identity (\ref{J2.5}) could be rewritten as follows
\cite{ZSSS}:
\begin{equation}
a\cdot (b\cdot (c\cdot d)) + (b\cdot (a\cdot c))\cdot d
 + c\cdot(b\cdot (a\cdot d))
 = (a\cdot b)\cdot (c\cdot d) + (a\cdot c)\cdot
(b\cdot d) + (c\cdot b) \cdot (a\cdot d).
                               \label{J2.6}
\end{equation}

Now, let  $C$ be a conformal algebra with $n$-products
$(\cdot \oo{n}\cdot )$,
$n\ge 0$,
and let
$*: C\otimes C \to H^{\otimes 2}\otimes _H C$
 be the pseudoproduct~(\ref{pseudoprod}).
Denote
$h_n = (-D)^{n}/n!$.
It is clear that $h_nh_m = \binom{n+m}{n}h_{n+m}$.

Theorems \ref{thm 3.1} and \ref{thm 3.2}
imply that $C$ is a Jordan conformal algebra
if and only if
$C$ satisfies the identities
\begin{equation}
a* b = (\sigma _{12}\otimes_H \idd_C)(b* a)
                               \label{J2.7}
\end{equation}
and
\begin{multline}
 (a * (b * (c * d)))   \\
 +  (\sigma _{12}\otimes _H \idd_C ) ((b * (a * c)) * d)
 + (\sigma _{13}\otimes _H \idd_C)(c* (b * (a * d)))      \\
=
((a * b)* (c* d ))
+  (\sigma _{23}\otimes _H \idd_C)((a* c)* (b * d))   \\
+ (\sigma _{13}\otimes _H \idd_C)((c * b)* (a * d)),
                                         \label{J2.8}
\end{multline}
where
$\sigma _{ij} = (i\,j)\in {\mathbb S}_4$.

Identity (\ref{J2.7}) is equivalent to the
conformal commutativity (\ref{ConfComm}),
so let us proceed with (\ref{J2.8}).
For example,
\begin{multline}
((b*(a*c))*d)
=
\bigg ( b*
\sum\limits_{n\ge 0} h_n\otimes 1 \otimes_H (a\oo{n} c)
 \bigg) *d
\\
=
\bigg(\sum\limits_{n,m\ge 0}
h_m\otimes h_n\otimes 1\otimes _H (b\oo{m}(a\oo{n} c))
\bigg)*d    \\
=
\sum\limits_{n,m,l\ge 0}
(h_m\otimes h_n \otimes 1\otimes 1)(\Delta^{(3)}(h_l)\otimes 1)
 \otimes_H (b\oo{m}(a\oo{n}c))\oo{l} d.
                                            \nonumber
\end{multline}
Hence,
\begin{multline}
(\sigma _{12}\otimes_H \idd_C)((b*(a*c))*d)
                                      \\
  =
\sum\limits_{\begin{array}{c}
           \scriptstyle  n,m,l\ge 0 \\
           \scriptstyle  s_1,s_2 \ge 0
           \end{array}}
\binom{n}{s_1}\binom{m}{s_2}
 (b\oo{m-s_2}(a \oo{n-s_1} c)\oo{l+s_1+s_2} d).
\end{multline}

By the same way,
one may proceed with other monomials in (\ref{J2.8})
and get the equivalent relation (more precisely,
this is a system of relations)
 in terms of conformal operations:
\begin{multline}
a\oo{n}(b\oo{m}(c\oo{l} d))
+
\sum\limits_{s_1,s_2 \ge 0}
 \binom{n}{s_1}\binom{m}{s_2}
 (b\oo{m-s_2}(a \oo{n-s_1} c)\oo{l+s_1+s_2} d) \\
+
c\oo{l} (b\oo{m}(a\oo{n} d))
=
\sum\limits_{s\ge 0}
 \binom{n}{s}
 (a\oo{n-s} b)\oo{m+s}(c\oo{l} d)  \\
+
\sum\limits_{s\ge 0}
 \binom{n}{s}
 (a\oo{n-s} c)\oo{l+s}(b\oo{m} d)
+
\sum\limits_{s\ge 0}
 \binom{l}{s}
 (c\oo{l-s} b)\oo{m+s} (a\oo{n} d)
                                            \label{J3.14}
 \end{multline}
for any $n,m,l\ge 0$.

\begin{prop}\label{prop 3.4}
A commutative conformal algebra $C$ is Jordan
if and only if
$C$ satisfies the identities~(\ref{J3.14}).
\qed
\end{prop}

\begin{defn}[see, e.g., \cite{Ku,S,ZSSS}]\label{defnAlter}
{\rm (i)}
 An algebra
 $(A, \cdot )$
is said to be
 {\em left\/}
 or
 {\em right alternative},
if for any
$a,b\in A$
we have
\begin{eqnarray}
 a^2\cdot b = a\cdot (a\cdot b),
                                \label{J4.1l}\\
 a\cdot b^2 = (a\cdot b)\cdot b.
                                \label{J4.1}
\end{eqnarray}

{\rm (ii)}
An anti-commutative algebra
$(A, \cdot)$
is said to be a Mal'cev algebra,
if it satisfies
the identity
\begin{equation}
 J(a,b,a\cdot c) = J(a,b,c)\cdot a,
                                       \label{J4.3}
\end{equation}
where
$J(a,b,c) = (a\cdot b)\cdot c - a\cdot (b\cdot c) + b\cdot (a\cdot c)$,
$a,b,c\in A$.
\end{defn}

The following statements are just corollaries of
Theorems \ref{thm 3.1} and \ref{thm 3.2}.

\begin{prop}\label{prop Alt}
A conformal algebra $C$ is left or right alternative
if and only if $C$ as a pseudoalgebra satisfies the relations
\begin{equation}
a*(b*c)-(a*b)*c
= (\sigma_{12}\otimes_H \idd_C)((b*a)*c - b*(a*c))
                                 \label{J4.2l}
\end{equation}
or
\begin{equation}
a*(b*c) - (a*b)*c
= (\sigma_{23}\otimes_H \idd_C)((a*c)*b - a*(c*b)),
                                \label{J4.2}
\end{equation}
respectively.
\end{prop}

\begin{prop}\label{prop Mal}
An anti-commutative conformal algebra $C$ is a
Mal'cev conformal algebra if and only if
it satisfies the identity
\begin{equation}
J^*(a,b,c*d) - J^*(a,b,c)*d
= (\sigma_{23}\otimes_H \idd_C)(a*J^*(c,b,d) - J^*(a*c, b,d)),
                                 \label{J4.Mal}
\end{equation}
where
$J^*(a,b,c)= (a*b)*c - a*(b*c) + (\sigma_{12}\otimes_H \idd_C)(b*(a*c))$.
\end{prop}

\begin{rem}
The ``pseudo"-Jacobian $J^*$ in (\ref{J4.Mal})
is mentioned to be defined with respect to Proposition \ref{prop 3.3}.
\end{rem}

The following statement describes
some elementary relations between
the considered varieties of pseudoalgebras.
To check these properties, one should perform
exactly the same computations as for
ordinary algebras, using Proposition~\ref{prop 3.3}.

\begin{prop}[c.f. \cite{BDK,K3}]\label{prop}
Let $H$ be a cocommutative Hopf algebra and let
$P$  be an $H$-pseudoalgebra with
a pseudoproduct~$*$.
Define the following $H$-bilinear maps on $P\otimes P$:
\[
[a * b]_{+} = a*b + (\sigma _{12}\otimes_H \idd_P)(b*a),
\quad
[a * b]_{-} = a*b - (\sigma _{12}\otimes_H \idd_P)(b*a).
\]
Denote by $P^{(\pm)}$
the same $H$-module $P$ endowed
with the pseudoproduct $[\cdot * \cdot]_{\pm})$.

{\rm (i)}
If $P$ is associative, then $P^{(-)}$ is Lie and $P^{(+)}$ is Jordan.

{\rm (ii)}
If $P$
is alternative, then
$P^{(-)}$ is Mal'cev and
$P^{(+)}$ is Jordan.
\end{prop}

It is also easy to write down the identities (\ref{J4.2l}),
(\ref{J4.2}) and (\ref{J4.Mal}) in terms of conformal products,
as it was done for Jordan identity.
But since the ``conformal''  form of identities is more
complicated than the ``pseudoalgebraic'' one (e.g., compare
(\ref{J2.8}) with (\ref{J3.14})), the language of pseudoproduct
seems to be more adequate even in the case of conformal algebras.

\begin{defn}
Let $H$ be a Hopf algebra.
An algebra $A$ (non-associative, in general)
endowed with homomorphism of algebras $\Delta_A : A\to H\otimes A$,
$\Delta_A(a)=a_{(1)}\otimes a_{(2)}$,
is said to be an $H$-{\em comodule algebra},
if
\begin{eqnarray}
&
  (\idd_H \otimes \Delta_A)\Delta_A(a) = (\Delta\otimes\idd_A)\Delta_A(a)
   = a_{(1)}\otimes a_{(2)}\otimes_{(3)}, \\
&
  \varepsilon(a_{(1)})a_{(2)} = a.
\end{eqnarray}
\end{defn}

The following statement shows how to construct conformal algebras
satisfying homogeneous multilinear identities.

\begin{prop}\label{prop 3.3x}
Let $H$ be a commutative and cocommutative Hopf algebra
and let $A$ be an $H$-comodule algebra.
Then the free $H$-module $P(A)=H\otimes A$
with the pseudoproduct~$*$ given by
\begin{equation}
 (h\otimes a)*(g\otimes b)
= hb_{(1)}\otimes ga_{(1)}\otimes_H (1\otimes a_{(2)}b_{(2)}),
\quad
h,g\in H,\ a,b\in A,
                                             \label{eq3.19}
\end{equation}
is an $H$-pseudoalgebra.
If $A$ satisfies an identity~(\ref{J3.8})
then the pseudoalgebra
$P(A)$ satisfies~(\ref{J3.8'}).
\end{prop}

\begin{proof}
Let $t(a_1,\dots ,a_n)$
be a non-associative word
obtained from $a_1\dots a_n$ by some bracketing.
It is sufficient to prove that
\begin{multline}
 t^*(a_1,\dots ,a_n)
                     \\
  =
\left (\bigoplus\limits_{k=1}^n
a_{1(k-1)}\dots a_{k-1(k-1)} a_{k+1(k)} \dots  a_{n(k)}
\right)
\otimes_H (1\otimes t(a_{1(n)},\dots ,a_{n(n)})).
                   \nonumber
\end{multline}
It could be easily done by induction on~$n$.
\end{proof}

In particular, let $H=\Bbbk[D]$, and let $A$ be an algebra.
Denote by
$A[t]$ the tensor product $A\otimes \Bbbk [t]$.
If $A$ satisfies a homogeneous multilinear identity
$f=0$ of type (\ref{J3.8}), then
$A[t]$ endowed with
\[
\Delta_{A[t]}: a\otimes t^{n} \to \sum\limits_{s\ge 0}
  \binom{n}{s} D^{s}\otimes (a\otimes t^{n-s})
\]
is an $H$-comodule algebra satisfying the same identity.
Hence, the pseudoalgebra
$P(A[t])$ constructed in Proposition~\ref{prop 3.3x}
satisfies the corresponding $f^*=0$ of type (\ref{J3.8'}).

There is a well-known fact (Artin's Theorem) in the theory
of alternative algebras. It states that any two elements
of an alternative algebra $A$ generate an associative
subalgebra of $A$. Let us show that this statement does not
hold for alternative conformal algebras.

Let $A$ be an alternative conformal algebra which is not associative.
Then $P(A[t])$ constructed above is an alternative
pseudoalgebra. Let us choose $a,b,c\in A $ such that
$\{a,b,c\}\equiv (ab)c - a(bc)\ne 0$, and fix $x,y\in P(A[t])$
as follows:
   $x=1\otimes (a\otimes t + b\otimes 1)$,
   $y=1\otimes (c\otimes 1)$.
Direct computation shows that
\begin{multline}
(x*y)*x - x*(y*x)
  \\
=(D\otimes 1\otimes 1)\otimes_H (1\otimes \{b,c,a\}\otimes 1)
+ (1\otimes 1\otimes D)\otimes_H (1\otimes\{a,c,b\}\otimes 1)\ne 0.
\nonumber
\end{multline}
Therefore, $x$ and $y$ do not generate an associative conformal subalgebra
of $P(A[t])$.

%\vfill\eject

\end{document}